\documentclass{article}
\usepackage{amssymb}

\input{tcilatex}
\begin{document}

\centerline{ Simplifying the axiomatization for ordered affine geometry
via a theorem prover}

\centerline{Dafa Li} \centerline{Dept of mathematical sciences}

\centerline{Tsinghua University, Beijing 100084, China}

\centerline{Abstract}

Jan von Plato proposed in 1998 an intuitionist axiomatization of ordered
affine geometry consisting of 22 axioms. It is shown that axiom I.7, which
is equivalent to a conjunction of four statements, two of which are
redundant, can be replaced with a simpler axiom, which is von Plato's
Theorem 3.10.


\section{Introduction}

Heyting proposed a constructive axiomatization for affine geometry \cite%
{Heyting}, by adopting the primitive notions of distinct points and distinct
lines. von Plato proposed 22 axioms for constructive ordered affine geometry 
\cite{jvp-98}, based on the five basic relations DiPt, DiLn, Undir, L-Apt,
and L-Con and four constructions ln($a,b$), pt($a,b$), par$(l,m)$, and rev($%
l $).

With the aid of the automatic theorem prover ANDP, we simplify in this paper
von Plato's axiomatization. Von Plato's axiom I.7 is shown to include some
redundancy and we show that it can be replaced with a simpler and more
intuitive new axiom, which is von Plato's Theorem 3.10 \cite{jvp-98}.

\section{Von Plato's axiomatization for ordered affine geometry}

In this paper, we want to replace axiom I.7 with a simpler and more
intuitive axiom. Only four axioms I.5-I.8 are of concern to the presented
work \cite{jvp-98}. Only a basic relation Undir$(l,m)$ and only a
construction rev$(l)$ appear in the axioms I.5-I.8. Undir$(l,m)$ means $l$
and $m$ are unequally directed lines. rev($l$) stands for the reverse of
line $l$. For readability, we list the four axioms as follows.

The axiom I.5: 
\begin{equation}
(\forall l)\thicksim Undir(l,l).
\end{equation}

The axiom I.6: 
\begin{equation}
(\forall l)(\forall m)(\forall n)[Undir(l,m)\rightarrow Undir(l,n)\vee
Undir(m,n)]  \label{I-6}
\end{equation}

The axiom I.7:%
\begin{eqnarray}
&&(\forall l)(\forall m)(\forall n)[Undir(l,m)\&Undir(l,rev(m))  \nonumber \\
&\rightarrow &Undir(l,n)\&Undir(l,rev(n))\vee Undir(m,n)\&Undir(m,rev(n))].
\label{I-7}
\end{eqnarray}

The axiom I.8:%
\begin{equation}
(\forall l)(\forall m)[Undir(l,m)\vee Undir(l,rev(m))]
\end{equation}

The convergence $Con(l,m)$ is defined by means of the basic relation Undir\
as $Undir(l,m)\&Undir(l,rev(m))$ \cite{jvp-98}. Thus, $Con(l,m)$ means that $%
l$ and $m$ are convergent lines. And then, axiom I.7 becomes 
\begin{equation}
(\forall l)(\forall m)(\forall n)[Con(l,m)\rightarrow Con(l,n)\vee Con(m,n)].
\label{conv}
\end{equation}

(\ref{conv}) describes the convergence of three lines (the direction is
ignored).

Axiom I.7 can be rewritten constructively equivalently as w1\&w2\&w3\&w4,
where wi, i=1,2,3,4, are listed in order as follows.

\begin{eqnarray}
Undir(l,m)\&Undir(l,rev(m)) &\rightarrow &Undir(l,n)\vee Undir(m,n),
\label{W-1} \\
Undir(l,m)\&Undir(l,rev(m)) &\rightarrow &Undir(l,n)\vee Undir(m,rev(n)),
\label{W-2} \\
Undir(l,m)\&Undir(l,rev(m)) &\rightarrow &Undir(l,rev(n))\vee Undir(m,n),
\label{W-3} \\
Undir(l,m)\&Undir(l,rev(m)) &\rightarrow &Undir(l,rev(n))\vee
Undir(m,rev(n)),  \label{W-4}
\end{eqnarray}%
where the prefix $(\forall l)(\forall m)(\forall n)$ for w1-w4 is omitted.

\section{Replacing axiom I.7 with a shorter formula}

By means of ANDP, after lots of trials, we chose SYM to replace axiom I.7,
where SYM\ stands for: 
\begin{equation}
(\forall l)(\forall m)[Undir(l,rev(m))\rightarrow Undir(m,rev(l))].
\end{equation}

SYM reads that for any lines $l$ and $m$, if $l$ and the reverse of $m$ are
unequally directed lines, then $m$ and the reverse of $l$ are unequally
directed lines. Here, SYM is referred as to the symmetry of the reverse
line. Clearly, SYM is shorter and more intuitive than axiom I.7. One notices
that there are two and respectively six occurrences of the predicate Undir
in SYM and in axiom I.7.

The proofs obtained via ANDP are classical. By means of the ANDP proofs, we
find the corresponding resolution proofs omitted here and corresponding
constructive proofs below.

\subsection{\textit{\ }Axioms I.5, I.6, and SYM constructively imply axiom
I.7}

To reduce the difficulty of deriving axiom I.7, we derive w1, w2, w3, and w4
from the axioms I.5, I.6, and SYM.

It is easily seen that w1 and w4 weakenings of I.6 and thus follow from I.6
alone.

\subsubsection{Proof of w2}

\textit{Lemma 1. }w2 can be constructively\ derived from the axioms I.5, I.6
and SYM.

Proof of (\ref{W-2}). Assume ($^{\ast }$) $Undir(l,m)$ and ($^{\ast \ast }$) 
$Undir(l,rev(m))$. Then, let us derive $Undir(l,n)\vee Undir(m,rev(n))$.

\begin{equation}
Undir(l,rev(m))\rightarrow Undir(l,n)\vee Undir(rev(m),n)\text{ \ \ \ \ \ \
\ \ (by I.6)}  \label{b1}
\end{equation}

\begin{equation}
Undir(l,n)\vee Undir(rev(m),n)\text{ \ \ \ \ (by (}^{\ast \ast }\text{) and (%
\ref{b1}))}  \label{b2}
\end{equation}

There are the following two cases for (\ref{b2}).

Case 1. $Undir(l,n)$. Therefore, $Undir(l,n)\vee Undir(m,rev(n))$.

Case 2. ($^{\ast \ast \ast }$) $Undir(rev(m),n)$. 
\begin{equation}
Undir(rev(m),n)\rightarrow Undir(rev(m),rev(m))\vee Undir(n,rev(m))\text{ \
\ \ \ (by I.6)}  \label{b4}
\end{equation}

\begin{equation}
Undir(rev(m),rev(m))\vee Undir(n,rev(m))\text{ \ \ \ \ \ \ \ (by (}^{\ast
\ast \ast }\text{) and (\ref{b4}))}  \label{b5}
\end{equation}

\begin{equation}
\sim Undir(rev(m),rev(m))\text{ \ \ \ \ \ \ \ \ (by I.5)}  \label{ad-1}
\end{equation}

\begin{equation}
Undir(n,rev(m))\text{ \ \ \ \ \ \ \ (by (\ref{b5}) and (\ref{ad-1})) }
\label{ad-2}
\end{equation}

\begin{equation}
Undir(n,rev(m))\rightarrow Undir(m,rev(n))\text{ \ \ \ \ \ \ (by SYM)}
\label{b6}
\end{equation}

\begin{equation}
Undir(m,rev(n))\text{ \ \ \ \ \ \ \ (by (\ref{ad-2}) and (\ref{b6}))}
\label{ad-3}
\end{equation}

\begin{equation}
Undir(l,n)\vee Undir(m,rev(n))\text{ \ \ \ \ \ \ (by (\ref{ad-3}))}
\end{equation}

Therefore, in all cases, we have $Undir(l,n)\vee Undir(m,rev(n))$. Q.E.D.

\subsubsection{Proof of w3}

\textit{Lemma 2}. w3 can be constructively\ derived from the axioms I.5, I.6
and SYM.

Proof of (\ref{W-3}). Assume ($^{\ast }$) $Undir(l,m)$ and ($^{\ast \ast }$) 
$Undir(l,rev(m))$. Then, let us derive $Undir(l,rev(n))\vee Undir(m,n)$.

\begin{equation}
Undir(l,rev(m))\rightarrow Undir(m,rev(l))\text{ \ \ \ \ \ \ (by SYM)}
\label{sy-1}
\end{equation}

\begin{equation}
Undir(m,rev(l))\text{ \ \ \ \ \ \ \ \ \ \ \ (by (}^{\ast \ast }\text{) and (%
\ref{sy-1}))}  \label{sy-2}
\end{equation}

\begin{equation}
Undir(m,rev(l))\rightarrow Undir(m,n)\vee Undir(rev(l),n)\text{ \ \ \ \ \ \
\ \ \ \ (by I.6)}  \label{c2}
\end{equation}

\begin{equation}
Undir(m,n)\vee Undir(rev(l),n)\text{ \ \ \ \ \ \ \ \ \ \ \ \ (by (\ref{sy-2}%
) and (\ref{c2}))}  \label{c3}
\end{equation}

There are the following two cases for (\ref{c3}).

Case 1. $Undir(m,n)$. Therefore, $Undir(l,rev(n))\vee Undir(m,n)$.

Case 2. ($^{\ast \ast \ast }$) $Undir(rev(l),n)$. 
\begin{equation}
Undir(rev(l),n)\rightarrow Undir(rev(l),rev(l))\vee Undir(n,rev(l))\text{ \
\ \ \ \ (by I.6)}  \label{c4}
\end{equation}

\begin{equation}
Undir(rev(l),rev(l))\vee Undir(n,rev(l))\text{ \ \ \ \ \ (by (}^{\ast \ast
\ast }\text{) and (\ref{c4})}  \label{c5}
\end{equation}

\begin{equation}
\sim Undir(rev(l),rev(l))\text{ \ \ \ \ \ \ \ \ \ \ \ (by I.5)}  \label{sy-3}
\end{equation}

\begin{equation}
Undir(n,rev(l))\text{ \ \ \ \ \ \ \ \ \ \ \ (by (\ref{sy-3}) and (\ref{c5}))}
\label{sy-4}
\end{equation}

\begin{equation}
Undir(n,rev(l))\rightarrow Undir(l,rev(n)\text{ \ \ \ \ \ \ \ (by SYM) }
\label{c6}
\end{equation}

\begin{equation}
Undir(l,rev(n)\text{ \ \ \ \ \ (by (\ref{sy-4}) and (\ref{c6}))}
\label{sy-5}
\end{equation}

\begin{equation}
Undir(l,rev(n))\vee Undir(m,n)\text{ \ \ \ \ \ \ \ (by (\ref{sy-5}))}
\end{equation}

Therefore, in all cases, we have $Undir(l,rev(n))\vee Undir(m,n)$. Q.E.D.

In light of Lemmas 1 and 2, we can conclude the following.

\textit{Theorem 1. }Axioms I.5, I.6, and SYM constructively imply axiom I.7.

\subsection{Conversely, the axioms I.5-I.8 constructively imply SYM}

\textit{Lemma 3}. SYM can be constructively derived from the axioms I.5,
I.6, I.8, and w2.

Proof. First, we want to indicate that SYM is von Plato's Theorem 3.10 \cite%
{jvp-98}. To derive SYM from the axioms I.5, I.6, I.8, and w2, we only need
to derive $Undir(m,rev(l))$ by assuming ($^{\ast }$) $Undir(l,rev(m))$. 
\begin{equation}
Undir(m,l)\vee Undir(m,rev(l))\text{ \ \ \ \ \ \ \ (by I.8)}  \label{sym-1}
\end{equation}

There are the following two cases for (\ref{sym-1}).

Case 1.$\ $($^{\ast \ast }$) $Undir(m,l)$. 
\begin{equation}
Undir(m,l)\rightarrow Undir(m,m)\vee Undir(l,m)\text{ \ \ \ \ \ \ (by I.6)}
\label{sym-2}
\end{equation}

\begin{equation}
Undir(m,m)\vee Undir(l,m)\text{ \ \ \ \ \ \ \ \ \ \ (by (}^{\ast \ast }\text{%
) and (\ref{sym-2}))}  \label{sym-3}
\end{equation}

\ 
\begin{equation}
\sim Undir(m,m)\text{ \ \ \ \ \ \ \ (by I.5)}  \label{sm-1}
\end{equation}

\begin{equation}
Undir(l,m)\text{ \ \ \ \ \ \ \ \ \ \ \ (by (\ref{sym-3}) and (\ref{sm-1}))}
\label{sm-2}
\end{equation}

\begin{equation}
Undir(l,m)\wedge Undir(l,rev(m))\rightarrow Undir(l,l)\vee Undir(m,rev(l))%
\text{ \ \ \ \ (by (\ref{W-2}))}  \label{sym-6}
\end{equation}

\begin{equation}
Undir(l,l)\vee Undir(m,rev(l))\text{ \ \ \ \ \ \ \ \ (by (}^{\ast }\text{), (%
\ref{sm-2}), and (\ref{sym-6}))}  \label{sym-7}
\end{equation}

\begin{equation}
\sim Undir(l,l)\text{ \ \ \ \ \ \ \ \ \ \ \ \ \ (by I.5)}  \label{sm-3}
\end{equation}

\begin{equation}
Undir(m,rev(l))\text{ \ \ \ \ \ \ \ \ \ \ \ \ (by (\ref{sym-7}) and (\ref%
{sm-3}))}
\end{equation}

Case 2. $Undir(m,rev(l))$. This is the desired result.

Therefore, in all cases, we have $Undir(m,rev(l))$. Q.E.D.

In light of Lemma 3, we can conclude the following.

\textit{Theorem 2. }The axioms I.5-I.8 constructively\ imply SYM.

In light of Theorems 1 and 2, we obtain the following.

\textit{Theorem 3. }The set of the axioms I.5, I.6, I.8 and SYM is
constructively equivalent to the set of the axioms I.5-I.8. In other words,
axiom I.7 is constructively equivalent to SYM, by using the axioms I.5, I.6,
and I.8. Thus, axiom I.7 can be replaced with SYM.

Acknowledgements

Thank the reviewer for the deep and helpful comments.

No funding; D. Li wrote the main manuscript text.; no conflict interest;

Data availability statement: All data generated or analysed during this
study are included in this published article.

\end{document}